# On Estimates of Exponential Sums

Nelson A. Carella, June 2002, NY

***Abstract:*** This paper introduces a general technique for estimating the absolute value of the exponential sum $S_k(s) = \Sigma_{0 \leq x < p} \exp(i2\pi s x^k / p)$ for all composite integers $k \neq q$, where $q$ is a prime power. The new estimate improves the classical estimate $\left| S_k(s) \right| \leq (k-1)p^{1/2}$ by a factor of about 2 or better (depending on $k$) for all sufficiently large primes $p$. Further, it will be shown that $\left| S_k(s) \right| \leq 2k^{3/4}p^{1/2}$ holds for a large class of integers $k$.

## 1 Introduction

Let $k, n, p \in \mathbb{N}$, $p$ prime, and $0 \neq s \in \mathbf{F}_p$. This paper introduces a general technique for estimating the absolute value of the exponential sum

$$S_k(s) = \sum_{x=0}^{p-1} e^{i2\pi s x^k / p} \qquad (1)$$

for all composite integers $k \neq q$, $q$ a prime power. The new estimate improves the classical estimate $\left| S_k(s) \right| \leq (k-1)p^{1/2}$ by a factor of about 2 or better (depending on $k$) for all sufficiently large primes $p$. Further, it will be shown that $\left| S_k(s) \right| \leq 2k^{3/4}p^{1/2}$ holds for a large class of integers $k$.

The new estimates are formally stated in Theorems 4 and 10. A few other known estimates are also stated here as references.

The classical estimate in Theorem 1 can be traced back to the 1800's. In the last decade of 1995 some improvements were accomplished by several authors. A notable case given in Theorem 2 improves the classical estimate by a factor of about $2^{1/2}$ uniformly for all even $k$, and all primes of the form $p = km + 1$, $m$ odd.

**Theorem 1.** If $p$ is a prime, then $\left|S_k(s)\right| \leq (k-1)p^{1/2}$ for all integer $k$. In addittion, if $k$ is odd, then $S_k \in \mathbb{R}$ is a real number.

**Theorem 2.** ([3]) Let $k$ be an even natural number. Suppose that $p$ is an odd prime number with $p \equiv 1 \bmod k$, $2k \nmid p-1$ and $p \nmid a$. Then

$$\left|S_k(s)\right| \leq 2^{-1/2}(k^2 - 2k + 2)^{1/2}p^{1/2}. \tag{2}$$

The most recent development on this topic is the following.

**Theorem 3.** Let $k \in \mathbb{N}$ be an integer, let $p$ be a prime, and $0 \neq s \in \mathbf{F}_p$. Then

$$\left|S_k(s)\right| \leq \begin{cases} kp^{1/2} & \text{if } 1 \leq k \leq p^{1/3}, \\ k^{5/8}p^{5/8} & \text{if } p^{1/3} \leq k \leq p^{1/2}, \\ k^{3/8}p^{3/4} & \text{if } p^{1/2} \leq k \leq p^{2/3}, \\ p & \text{if } p^{2/3} \leq k \leq p. \end{cases} \tag{3}$$

This estimate is derived from the estimated number of solutions of the congruence $x^k + y^k \equiv n \bmod p$. The present form appears in [1], see also [4] for a generalization.

The contribution here is the following.

**Theorem 4.** If $k \in \mathbb{N}$ is a divisor of $p-1$, $p$ a prime, then

(i) $\left|S_k(s)\right| \leq \left(d - 1 + (n-1)\sqrt{d(1 + dnp^{-1/2})}\right)p^{1/2}$, if $k = dn$, and $\gcd(d, n) = 1$.

(ii) $\left|S_k(s)\right| \leq \left(d - 1 + \sqrt{d(1 + dp^{-1/2})}\right)p^{1/2}$, if $k = 2d$, and $\gcd(d, 2) = 1$.

(iii) $\left|S_k(s)\right| \leq \left(d - 1 + 3\sqrt{d(1 + 2^{1/2}\,p^{-1/2})}\right)p^{1/2}$, if $k = 4d$, and $\gcd(d, 4) = 1$.

Write $\left|S_k(s)\right| \leq c(k)p^{1/2}$. The classical estimate has the coefficient $c_0(k) = k - 1$, this is nontrivial for all $k < p^{1/2}$. The estimate reported in Theorem 2 has the coefficient $c_1(k) = 2^{-1/2}(k^2 - 2k + 2)^{1/2}$ for all primes of the form $p = km + 1$, $m$ odd. The new estimate in Theorem 4 has a coefficient of the form

$$c_d(k) = d - 1 + (n-1)(d(1 + kp^{-1/2}))^{1/2} \tag{4}$$

for fixed $k = dn$, and $p$ prime, (there are certain conditions on $k$ and the divisor $d$). Since the ratio $c_0(k)/c_1(k) \approx 2^{1/2}$, the coefficient $c_1(k)$ improves the classical coefficient $c_0(k)$ by a factor of about $2^{1/2}$. Likewise, the ratio $c_0(k)/c_d(k) \approx (d^{1/2}n)/(d^{1/2} + n)$ implies that the coefficient $c_d(k)$ improves the classical coefficient $c_0(k)$ by a factor of about $(d^{1/2}n)/(d^{1/2} + n)$. This assumes that $d < k < p^{1/2}$, and $p$ is moderately large so that $(1 + kp^{-1/2})^{1/2} \approx 1$.

The next section deals with auxiliary results, and an estimate of the nonlinear exponential sum $T_d(\chi,s)$. And the last section covers the proof of Theorem 4, and another result for nonprimelike integers $k$, (integers with two or more large relatively prime factors).

## 2 Some Properties of $S_k(s)$ and $T_d(\chi,s)$

The integer $k$ is restricted to divisors of $p - 1$ due to the fact that $\gcd(k, p - 1) = d \geq 1$ implies that $S_k(s) = S_d(s)$.

***Theorem* 5.** The second moment of $S_k(s)$ is $\sum_{0 \leq s < p} |S_k(s)|^2 = kp(p - 1)$.

This indicates that the real numbers $|S_k(s)|$, $0 \neq s$, are clustered around the value $(kp)^{1/2}$.

The *character spectrum* of the exponential sum $S_k(s)$ is a decomposition of $S_k(s)$ as a sum of closely related exponential sums. Under this scheme $S_k(s)$ is viewed as a function on $\mathbf{F}_p$, and the spectrum is a sort of transformation of $S_k(s)$. The new variable is $d \mid k$.

***Definition* 6.** Let $k = dn \geq 1$, and let $p$ be a prime. If $\chi$ is a multiplicative characters of order $\mathrm{ord}(\chi) = n$, then the character spectrum of $S_k(s)$ is defined by

$$\sum_{x=0}^{p-1} e^{i2\pi sx^k / p} = \sum_{e=0}^{n-1} \sum_{x=0}^{p-1} \chi^e(x) e^{i2\pi sx^d / p} \ . \tag{5}$$

The character spectrum is given in terms of the exponential sum $T_d(\chi^e,s) = \sum_{0 \leq x < p} \chi^e(x)\exp(i2\pi sx^d/p)$. The next result will be used to estimate the value $| T_d(\chi^e,s) |$.

***Lemma* 7** Let $k = dn$, $p$ a prime, and let $\chi \neq 1$ be a multiplicative character of order $n$. If $\gcd(d, n) = 1$, then

$$\left| \sum_{z^d \neq 0, 1} \chi^e(z^{-1}) \right| \leq dn \ , \qquad e = 1, 2, \ldots, n - 1. \tag{6}$$

Proof: Assume $k = dn$ is a divisor of $p - 1$, and put $N = (p - 1) / d$. Now write the nonzero elements of $\mathbf{F}_p$ as $z = \tau^{sd + tN}$, where $\tau$ is a generator of $\mathbf{F}_p$, and $0 \leq s < N$, $0 \leq t < d$. The hypothesis $\gcd(d, n) = 1$ implies that $2de / n \neq$ even integer, $0 < e < n$, so $\chi^e$

$\neq 1$. Moreover, $z^d \neq 0, 1$ implies that $s \neq 0$. Accordingly the sum becomes

$$\sum_{z^d \neq 0,1} \chi^e(z^{-1}) = d \sum_{s=1}^{N-1} e^{-i2\pi s de/n} ,$$ (7)

an incomplete nontrivial character sum. To complete the analysis put $N - 1 = an + b$, $0 \leq b < n$, and simplify to obtain

$$\sum_{z^d \neq 0,1} \chi^e(z^{-1}) = d \sum_{s=1}^{b} \overline{\chi}^e(\tau^s) .$$ (8)

The claim follows from this. ∎

Since $p$ is a prime this is not the sharpest estimate, but it is very useful. In many cases the sum can be evaluated exactly or estimated very closely. For example, for $k = 2d$, and $k = 4d$, ($d$ odd), the estimates are

$$\left| \sum_{z^d \neq 0,1} \chi(z^{-1}) \right| = d \ \text{ and } \ \left| \sum_{z^d \neq 0,1} \chi(z^{-1}) \right| \leq 2^{1/2} d .$$ (9)

**Theorem 8.** The second moment of $| T_d(\chi,s) |$ is $\sum_{0 \leq s < p} | T_d(\chi,s) |^2 = dp(p-1)$.

Proof: By definition, and the change of variable $y = xz$, the second moment is precisely

$$\sum_{s=0}^{p-1} \left| T_k(\chi,s) \right|^2 = \sum_{s=0}^{p-1} \sum_{y=1}^{p-1} \overline{\chi}(y) e^{-i2\pi s y^d / p} \sum_{x=1}^{p-1} \chi(x) e^{i2\pi s x^d / p}$$

$$= \sum_{z=1}^{p-1} \sum_{x=1}^{p-1} \chi(z^{-1}) \sum_{s=0}^{p-1} e^{i2\pi s x^d (1 - z^d) / p} .$$ (10)

Since only the roots of $1 - z^d = 0$ contribute to the total, the expression reduces to

$$\sum_{s=0}^{p-1} \left| T_d(\chi,s) \right|^2 = \left( 1 + \chi(\omega) + \chi(\omega^2) + \cdots + \chi(\omega^{d-1}) \right) p(p-1) = dp(p-1),$$ (11)

where $\omega$ is a primitive $d$th root of unity, and $\chi(\omega^i) = 1$, $0 \leq i < d$. ∎

The average value of the sum $| T_d(\chi,s) |$ over $s \neq 0$ is $(dp)^{1/2}$. It seems that the absolute value of any sum $T_d(\chi,s)$ with the parameters $d$, and $n$ such that $\gcd(d, n) \neq 1$ can exceed the average value $(dp)^{1/2}$, at least for some primes $p$ and/or some $s \neq 0$.

However, restricting the parameters $d$, and $n$ to $\gcd(d, n) = 1$, one can obtain an estimate of $|T_d(\chi,s)|$ very close to the average value.

Let $R = \{\, x^d : 0 \neq x \in \mathbf{F}_p \,\}$, $R_1 = gR_0$, $R_2 = g^2 R_0$, ..., $R_{d-1} = g^{d-1} R_0$, $g$ primitve modulo $p$. It is clear that the map $s \rightarrow |\,T_d(\chi,s)\,|$ is constant on any of the $d$th power residues classes $R_j = g^j R_0$. Specifically, $|\,T_d(\chi)\,| = |\,T_d(\chi,1)\,| = |\,T_d(\chi,s)\,|$ for any $d$th power residue $s \in R$ since $s = x^d$ for some $x \neq 0$.

**Theorem 9.** Let $d \in \mathbb{N}$ be a divisor of $p - 1$, $p$ prime, and let $\chi \neq 1$ be a multiplicative characters of order $\mathrm{ord}(\chi) = n$. Then

( 1 ) $\left|T_d(\chi^e,s)\right|^2 \leq d(1 + dn p^{-1/2})\, p$ , if $k = dn$, and $\gcd(d, n) = 1$, $0 < e < n$.

( 2 ) $\left|T_d(\chi^e,s)\right|^2 \leq d(1 + d p^{-1/2})\, p$ , if $k = 2d$, and $\gcd(d, 2) = 1$, $0 < e < n$.

( 3 ) $\left|T_d(\chi^e,s)\right|^2 \leq d(1 + 2^{1/2} p^{-1/2})\, p$ , if $k = 4d$, and $\gcd(d, 4) = 1$, $0 < e < n$.

Proof: Let $d > 1$, and fix a coset $R_j$, some $0 \leq j < n$. The sum of the square absolute values $|T_d(\chi,s)|^2$ over the coset $R_j$ is computed using two slightly different methods. The first summation method yields

$$\sum_{s \in R_j} \left|T_d(\chi,s)\right|^2 = \left(\frac{p-1}{d}\right)\left|T_d(\chi,s)\right|^2 . \qquad (12)$$

The second summation method yields

$$\sum_{s \in R_j} \left|T_d(\chi,s)\right|^2 = \sum_{s \in R_j}\sum_{z=1}^{p-1}\chi(z^{-1})\sum_{x=1}^{p-1} e^{i2\pi s x^d (1 - z^d)/p} . \qquad (13)$$

Removing the contribution from the roots of $1 - z^d = 0$ on the previous line and combining the two equations return

$$\left(\frac{p-1}{d}\right)\left|T_d(\chi,s)\right|^2 = (p-1)^2 + \sum_{s \in R_j,\, z^d \neq 0,1}\sum_{z=1}^{p-1}\chi(z^{-1})\sum_{x=1}^{p-1} e^{i2\pi s x^d (1 - z^d)/p} . \qquad (14)$$

The hypothesis $\gcd(d, n) = 1$ implies that $\chi^e \neq 1$, since $2de / n \neq$ even integer, $0 < e < n$, (this is required to use Lemma 7). Thus applying the triangle inequality and substituting the appropriate estimates one obtains

$$\left| T_d(\chi, s) \right|^2 = d(p-1) + \frac{d}{p-1} \sum_{s \in R_j} \sum_{z^d \neq 0,1} \chi(z^{-1}) \sum_{x=1}^{p-1} e^{i2\pi s x^d (1-z^d)/p}$$

$$\leq d(p-1) + \left( \frac{d}{p-1} \right) \left( \frac{p-1}{d} \right) (dn)((d-1)p^{1/2}) \qquad (15)$$

$$\leq d(1 + dnp^{-1/2})p.$$

The proofs of the last two estimates (2) and (3) use the exact value or a sharp estimate of the incomplete character sum. ∎

Note that in the case $d = 1$, the penultimate equation reduces to the standard Gaussian result $\left| T_1(\chi, s) \right|^2 = \left| G_n(\chi, s) \right|^2 = p$.

## 3 The Main Results

The next few lines completes Theorem 4.

Proof: From the spectrum decomposition of $S_k(s)$, and Theorem 9, it follows that

$$\left| S_k(s) \right| \leq \left| T_d(1, s) \right| + \left| T_d(\chi, s) \right| + \left| T_d(\chi^2, s) \right| + \cdots + \left| T_d(\chi^{n-1}, s) \right|$$

$$\leq [d - 1 + (n-1)(d(1 + kp^{-1/2}))^{1/2}]p^{1/2}. \qquad ∎$$

Since the function $f_k(x) = x - 1 + kx^{-1/2} - x^{1/2}$ has a minimal for some $x_0$ in the interval $(k^{1/2}, k^{3/4})$, using this method, the best estimates of $\left| S_k(s) \right|$ are realized with nonprimelike integers $k$, (integers with two or more large relatively prime factors). Under this condition there exists one or more factors $d$ of $k$ close to $x_0$, so the corresponding estimate $\left| S_k(s) \right|$ is close to the minimal $f_k(x_0)$ of $f_k(x)$. In contrast, the worst estimates of $\left| S_k(s) \right|$ are realized with prime power integers $k$, and primelike integers, for exanple, $k = 2n$, $n$ prime.

***Theorem 10.*** Suppose that the integer $k \in (17, p^{1/2})$ is not a prime power, and that it has a factor $d$ in the range $(k^{1/2}, k^{3/4})$, with $k = dn$, and $\gcd(d, n) = 1$. Then

$$\left| S_k(s) \right| \leq 2k^{3/4}p^{1/2}. \qquad (16)$$

Proof: Let $d = k^{1/2+\delta}$, $0 \leq \delta \leq 1/4$, and let $\xi = (1 + kp^{-1/2})^{1/2} \leq 2^{1/2}$. Then by Theorem 4, one has

$$\left| S_k(s) \right| \leq [k^{1/2+\delta} - 1 + (k^{1/2-\delta} - 1)(k^{1/2+\delta}(1 + kp^{-1/2}))^{1/2}]p^{1/2}$$

$$\leq k^{1/2}(k^\delta + \xi k^{1/4-\delta/2})p^{1/2}.$$

Since $k^\delta + \xi k^{1/4-\delta/2} = k^{1/4}(k^{\delta-1/4} + \xi k^{-\delta/2}) < 2k^{1/4}$, and $k > 17$, the claim follows. ∎

This last expression has some resemblance to the conjectured estimate

$$\left| S_k(s) \right| \leq (kp)^{1/2+\varepsilon}, \ \varepsilon > 0, \tag{17}$$

described in [3].